\documentclass[10pt]{article}
\usepackage{cite}
\usepackage{mathrsfs}
\usepackage{amsfonts}
\usepackage{amsmath}
\usepackage{amsfonts,amssymb,color}
\usepackage{dsfont}
\usepackage{curves}
\usepackage{mathrsfs}
\usepackage{pifont}
\usepackage{amssymb}
\allowdisplaybreaks

\numberwithin{equation}{section}

\date{}

\def\BigRoman{\uppercase\expandafter{\romannumeral\number\count 255 }}
\def\Romannumeral{\afterassignment\BigRoman\count255=}

\begin{document}
\title{Sufficient conditions for even factors in graphs
}
\author{\small Sizhong Zhou\footnote{Corresponding author. E-mail address: zsz\_cumt@163.com (S. Zhou)}, Qiuxiang Bian, Jiancheng Wu\\
\small  School of Science, Jiangsu University of Science and Technology,\\
\small  Zhenjiang, Jiangsu 212100, China\\
}

\maketitle
\begin{abstract}
\noindent Let $G$ be a graph. We denote by $e(G)$ and $\rho(G)$ the size and the spectral radius of $G$. A spanning subgraph $F$ of $G$ is called an even factor of $G$ if $d_F(v)\in\{2,4,6,\ldots\}$
for every $v\in V(G)$. Yan and Kano provided a sufficient condition using the number of odd components in $G-S$ for a graph $G$ of even order to contain an even factor, where $S$ is a vertex subset
of $G$ [Z. Yan, M. Kano, Strong Tutte type conditions and factors of graphs, Discuss. Math. Graph Theory 40 (2020) 1057--1065]. In this paper, motivated by Yan and Kano's above result, we present
some tight sufficient conditions to guarantee that a connected graph $G$ with the minimum degree $\delta$ contains an even factor with respect
to its size and spectral radius.
\\
\begin{flushleft}

(2020) Mathematics Subject Classification: 05C50, 05C70
\end{flushleft}
\end{abstract}

\section{Introduction}

In this paper, we deal with finite undirected graphs without loops and multiple edges. Let $G$ be a graph with vertex set $V(G)$ and edge set $E(G)$. We denote by $|V(G)|=n$ and $|E(G)|=e(G)$
the order and size of $G$, respectively. For $v\in V(G)$, let $d_G(v)$ denote the degree of $v$ in $G$. Let $\delta(G)=\min\{d_G(v):v\in V(G)\}$. The number of odd components in $G$ is denoted
by $o(G)$. For a vertex subset $S$ of $G$, let $G[S]$ and $G-S$ denote the subgraphs of $G$ induced by $S$ and $V(G)\setminus S$, respectively. The complete graph of order $n$ is denoted by
$K_n$. For two vertex-disjoint graphs $G_1$ and $G_2$, the union of $G_1$ and $G_2$ is denoted by $G_1\cup G_2$, and the join $G_1\vee G_2$ is the graph obtained from $G_1\cup G_2$ by adding
all possible edges between $V(G_1)$ and $V(G_2)$.

Let $A(G)$ denote the adjacency matrix of $G$. The largest eigenvalue of $A(G)$ is called the spectral radius of $G$, and denoted by $\rho(G)$. The study on the spectral radius of a graph has
attracted much attention \cite{GA,B,F,Wc,WZ,Zs,ZW,Zt,ZZL,ZZS,ZSL,WZL} in recent years.

A spanning subgraph $F$ of $G$ is called an even factor of $G$ if $d_F(v)\in\{2,4,6,\ldots\}$ for every $v\in V(G)$. In particular, an even factor $F$ is called a 2-factor if $d_F(v)=2$ for all
$v\in V(G)$.

Tutte \cite{T} presented a characterization for a graph containing a 2-factor. Ota and Tokuda \cite{OT} proved that a $K_{1,n}$-free graph $G$ with $\delta(G)\geq2n-2$ has a 2-factor, where
$n\geq3$ is an integer. Ryjáček, Saito and Schelp \cite{RSS} claimed that a $K_{1,3}$-free graph $G$ contains a 2-factor if and only if the closure of $G$ contains a 2-factor. Chen and Chen
\cite{CC} obtained a sufficient condition for a connected graph of order $n\geq3$ to contain a 2-factor. Fan and Lin \cite{FL} showed a spectral radius condition for the existence of 2-factors
in 1-binding graphs. Xiong \cite{X} proved two necessary and sufficient conditions for a special graph to possess an even factor. Steffen and Wolf \cite{SW} claimed that every $k$-critical graph
with at most $2k-6$ vertices of degree 2 contains an even factor. Fleischner \cite{F1} proved that a 2-edge-connected graph $G$ with $\delta(G)\geq3$ contains an even factor. Lv and Xiong \cite{LX}
posed a characterization of a graph $G$ such that the $n$-times iterated line graph $L^{n}(G)$ has an even factor with at most $k$ components. Moreover, they use this above result to establish
some upper bounds for the minimum number of components of even factors in $L^{n}(G)$. Kobayashi and Takazawa \cite{KT}, Zhang and Xiong \cite{ZX} provided some results on the existence of even
factors in graphs. Yan and Kano \cite{YK} provided sufficient conditions using the number of odd components in $G-S$ for a graph $G$ to contain an even factor, where $S\subseteq V(G)$. We refer
the reader to \cite{HZ,ZWH,ZZ,ZXS,Zr,Wr,Ws} for more results on graph factors.

Motivated by \cite{YK}, in this paper, we investigate some tight sufficient conditions to guarantee that a connected graph $G$ with the minimum degree $\delta$ contains an even factor with respect
to its size and spectral radius, and obtain the following two theorems.

\medskip

\noindent{\textbf{Theorem 1.1.}} Let $G$ be a connected graph of even order $n\geq\max\{6\delta-4,\frac{1}{6}(\delta^{2}+7\delta+4)\}$ with minimum degree $\delta\geq2$. If
$$
e(G)\geq e(K_{\delta}\vee(K_{n-2\delta+1}\cup(\delta-1)K_1)),
$$
then $G$ contains an even factor, unless $G=K_{\delta}\vee(K_{n-2\delta+1}\cup(\delta-1)K_1)$.

\medskip

\noindent{\textbf{Theorem 1.2.}} Let $G$ be a connected graph of even order $n\geq\max\{5\delta-3,\frac{1}{3}\delta^{2}+\delta\}$ with minimum degree $\delta\geq2$. If
$$
\rho(G)\geq\rho(K_{\delta}\vee(K_{n-2\delta+1}\cup(\delta-1)K_1)),
$$
then $G$ contains an even factor, unless $G=K_{\delta}\vee(K_{n-2\delta+1}\cup(\delta-1)K_1)$.

\section{Preliminary lemmas}

In this section, we introduce some preliminary lemmas, which are very useful in the proofs of our main theorems. Yan and Kano \cite{YK} showed a sufficient condition for the existence of even factors in
graphs.

\medskip

\noindent{\textbf{Lemma 2.1}} (Yan and Kano \cite{YK}). Let $G$ be a graph of even order $n$. Then $G$ contains an even factor if
$$
o(G-S)<|S|
$$
for all $S\subseteq V(G)$ with $|S|\geq2$.

\medskip

\noindent{\textbf{Lemma 2.2}} (Zheng, Li, Luo and Wang \cite{ZLLW}). Let $\sum\limits_{i=1}^{t}n_i=n-s$. If $n_1\geq n_2\geq\cdots\geq n_t\geq p\geq1$ and $n_1<n-s-p(t-1)$, then
$$
e(K_s\vee(K_{n_1}\cup K_{n_2}\cup\cdots\cup K_{n_t}))<e(K_s\vee(K_{n-s-p(t-1)}\cup(t-1)K_p)).
$$

\medskip

\noindent{\textbf{Lemma 2.3}} (Li and Feng \cite{LF}). Let $H$ be a subgraph of a connected graph $G$. Then
$$
\rho(G)\geq\rho(H),
$$
with equality occurring if and only if $G=H$.

\medskip

\noindent{\textbf{Lemma 2.4}} (Miao, Li and Wei \cite{MLW}). Let $n_1\geq n_2\geq\cdots\geq n_t\geq p$ be $(t+1)$ positive
integers with $\sum\limits_{i=1}^{t}n_i=n-s$ and $n_1<n-s-p(t-1)$. Then
$$
\rho(K_s\vee(K_{n_1}\cup K_{n_2}\cup\cdots\cup K_{n_t}))<\rho(K_s\vee(K_{n-s-p(t-1)}\cup(t-1)K_p)).
$$

\medskip

For a real $n\times n$ matrix 
\begin{align*}
M=\left(
  \begin{array}{cccc}
    M_{11} & M_{12} & \cdots & M_{1r}\\
    M_{21} & M_{22} & \cdots & M_{2r}\\
    \vdots & \vdots & \ddots & \vdots\\
    M_{r1} & M_{r2} & \cdots & M_{rr}\\
  \end{array}
\right),
\end{align*}
its rows and columns are partitioned by a partition $\pi:\mathcal{N}=\mathcal{N}_1\cup\mathcal{N}_2\cup\cdots\cup\mathcal{N}_r$, where $\mathcal{N}=\{1,2,\ldots,n\}$. Let $m_{ij}$ denote the average row
sum of the block $M_{ij}$ of $M$ for any $1\leq i,j\leq r$. Then the quotient matrix of $M$ is defined by $M_{\pi}=(m_{ij})_{r\times r}$. The partition $\pi$ is equitable if every block $M_{ij}$ of $M$
has constant row sum $m_{ij}$ for $1\leq i,j\leq r$.

\medskip

\noindent{\textbf{Lemma 2.5}} (You, Yang, So and Xi \cite{YYSX}). Let $M$ be a real $n\times n$ matrix with an equitable partition $\pi$, and let $M_{\pi}$ be the corresponding quotient matrix. Then each
eigenvalue of $M_{\pi}$ is an eigenvalue of $M$. Furthermore, if $M$ is nonnegative and irreducible, then the largest eigenvalues of $M$ and $M_{\pi}$ are equal.

\medskip

\section{The proof of Theorem 1.1}

\noindent{\it Proof of Theorem 1.1.} Suppose, to the contrary, that $G$ contains no even factor. In terms of Lemma 2.1, there exists some subset $S\subseteq V(G)$ such that $s=|S|\geq2$ and $o(G-S)\geq s$.
Obviously, $G$ is a spanning subgraph of $G_1=K_s\vee(K_{n_1}\cup K_{n_2}\cup\cdots\cup K_{n_s})$ for some odd integers $n_1\geq n_2\geq\cdots\geq n_s\geq1$ with $\sum\limits_{i=1}^{s}n_i=n-s$. Thus, we
conclude
\begin{align}\label{eq:3.1}
e(G)\leq e(G_1),
\end{align}
with equality occurring if and only if $G=G_1$. The following proof will be divided into three cases according to the value of $s$.

\noindent{\bf Case 1.} $s\geq\delta+1$.

Let $G_2=K_s\vee(K_{n-2s+1}\cup(s-1)K_1)$, where $n\geq2s$. By virtue of Lemma 2.2, we possess
\begin{align}\label{eq:3.2}
e(G_1)\leq e(G_2),
\end{align}
where the equality occurs if and only if $(n_1,n_2,\ldots,n_s)=(n-2s+1,1,\ldots,1)$.

Let $G_*=K_{\delta}\vee(K_{n-2\delta+1}\cup(\delta-1)K_1)$. A simple computation yields that
\begin{align*}
e(G_*)-e(G_2)=&\binom{n-\delta+1}{2}+\delta(\delta-1)-\binom{n-s+1}{2}-s(s-1)\\
=&\frac{1}{2}(s-\delta)(2n-3s-3\delta+3)\\
=&\frac{1}{4}(s-\delta)(4n-6s-6\delta+6)\\
\geq&\frac{1}{4}(s-\delta)(n-6\delta+6) \ \ \ \ \ (\mbox{since} \ s\geq\delta+1 \ \mbox{and} \ n\geq2s)\\
>&0 \ \ \ \ \ (\mbox{since} \ s\geq\delta+1 \ \mbox{and} \ n\geq6\delta-4).
\end{align*}
This implies that $e(G_2)<e(G_*)$. Together with \eqref{eq:3.1} and \eqref{eq:3.2}, we deduce
$$
e(G)\leq e(G_1)\leq e(G_2)<e(G_*)=e(K_{\delta}\vee(K_{n-2\delta+1}\cup(\delta-1)K_1)),
$$
which contradicts $e(G)\geq e(K_{\delta}\vee(K_{n-2\delta+1}\cup(\delta-1)K_1))$.

\noindent{\bf Case 2.} $s=\delta$.

Recall that $G_1=K_s\vee(K_{n_1}\cup K_{n_2}\cup\cdots\cup K_{n_s})$. In terms of Lemma 2.2 and $s=\delta$, we infer
\begin{align}\label{eq:3.3}
e(G_1)\leq e(K_{\delta}\vee(K_{n-2\delta+1}\cup(\delta-1)K_1)),
\end{align}
where the equality holds if and only if $G_1=K_{\delta}\vee(K_{n-2\delta+1}\cup(\delta-1)K_1)$. It follows from \eqref{eq:3.1} and \eqref{eq:3.3} that
$$
e(G)\leq e(K_{\delta}\vee(K_{n-2\delta+1}\cup(\delta-1)K_1)),
$$
with equality if and only if $G=K_{\delta}\vee(K_{n-2\delta+1}\cup(\delta-1)K_1)$, a contradiction.

\noindent{\bf Case 3.} $s\leq\delta-1$.

Let $G_3=K_s\vee(K_{n-s-(\delta+1-s)(s-1)}\cup(s-1)K_{\delta+1-s})$. Recall that $G$ is a spanning subgraph of $G_1=K_s\vee(K_{n_1}\cup K_{n_2}\cup\cdots\cup K_{n_s})$ for some odd integers
$n_1\geq n_2\geq\cdots\geq n_s$ with $\sum\limits_{i=1}^{s}n_i=n-s$. We easily see $n_s\geq\delta+1-s$ because $\delta(G_1)\geq\delta(G)=\delta$. Using Lemma 2.2, we possess
\begin{align}\label{eq:3.4}
e(G_1)\leq e(G_3),
\end{align}
with equality occurring if and only if $(n_1,n_2,\ldots,n_s)=(n-s-(\delta+1-s)(s-1),\delta+1-s,\ldots,\delta+1-s)$.

Recall that $G_*=K_{\delta}\vee(K_{n-2\delta+1}\cup(\delta-1)K_1)$. We easily get
$$
e(G_3)=\binom{n-(\delta+1-s)(s-1)}{2}+s(s-1)(\delta+1-s)+(s-1)\binom{\delta+1-s}{2}
$$
and
$$
e(G_*)=\binom{n-\delta+1}{2}+\delta(\delta-1).
$$
By a simple computation, we obtain
\begin{align}\label{eq:3.5}
e(G_*)-e(G_3)=&\frac{1}{2}((n-\delta)(n-\delta+1)+2\delta(\delta-1)\nonumber\\
&-(n-(\delta+1-s)(s-1))(n-(\delta+1-s)(s-1)-1)\nonumber\\
&-2s(s-1)(\delta+1-s)-(s-1)(\delta+1-s)(\delta-s))\nonumber\\
=&\frac{1}{2}(\delta-s)(s^{3}-(\delta+5)s^{2}+(2n+\delta+7)s-4n+3\delta-3).
\end{align}
Let $f(x)=x^{3}-(\delta+5)x^{2}+(2n+\delta+7)x-4n+3\delta-3$. Then
$$
f'(x)=3x^{2}-2(\delta+5)x+2n+\delta+7.
$$
Notice that the symmetry axis of $f'(x)$ is $x=\frac{\delta+5}{3}$. Hence, we infer $f'(x)\geq f'(\frac{\delta+5}{3})=\frac{1}{3}(6n-\delta^{2}-7\delta-4)\geq0$ by $n\geq\frac{1}{6}(\delta^{2}+7\delta+4)$,
which implies that $f(x)$ is increasing for $x\in[3,\delta-1]$. Combining this with $n\geq6\delta-4$ and $s\in[3,\delta-1]$, we have
\begin{align}\label{eq:3.6}
f(s)\geq f(3)=2n-3\delta\geq2(6\delta-4)-3\delta=9\delta-8>0.
\end{align}

It follows from \eqref{eq:3.5} and \eqref{eq:3.6} that
\begin{align}\label{eq:3.7}
e(G_3)<e(G_*)
\end{align}
for $s\in[3,\delta-1]$. We check that the inequality \eqref{eq:3.7} also holds for $s=2$. In view of \eqref{eq:3.1}, \eqref{eq:3.4} and \eqref{eq:3.7}, we deduce
$$
e(G)\leq e(G_1)\leq e(G_3)<e(G_*)=e(K_{\delta}\vee(K_{n-2\delta+1}\cup(\delta-1)K_1)),
$$
which contradicts $e(G)\geq e(K_{\delta}\vee(K_{n-2\delta+1}\cup(\delta-1)K_1))$. Theorem 1.1 is verified. \hfill $\Box$

\section{The proof of Theorem 1.2}

\noindent{\it Proof of Theorem 1.2.} Suppose, to the contrary, that $G$ contains no even factor. Then using Lemma 2.1, we conclude $o(G-S)\geq|S|$ for some subset $S\subseteq V(G)$ with $|S|\geq2$.
Let $|S|=s$. Then $G$ is a spanning subgraph of $G_1=K_s\vee(K_{n_1}\cup K_{n_2}\cup\cdots\cup K_{n_s})$, where $n_1\geq n_2\geq\cdots\geq n_s\geq1$ are odd integers and $\sum\limits_{i=1}^{s}n_i=n-s$.
By Lemma 2.3, we possess
\begin{align}\label{eq:4.1}
\rho(G)\leq\rho(G_1),
\end{align}
where the equality occurs if and only if $G=G_1$. We proceed by discussing the following three cases.

\noindent{\bf Case 1.} $s\geq\delta+1$.

Let $G_2=K_s\vee(K_{n-2s+1}\cup(s-1)K_1)$, where $n\geq2s$. In view of Lemma 2.4, we get
\begin{align}\label{eq:4.2}
\rho(G_1)\leq\rho(G_2),
\end{align}
where the equality holds if and only if $G_1=G_2$. For the partition $V(G_2)=V(K_s)\cup V(K_{n-2s+1})\cup V((s-1)K_1)$, the quotient matrix of $A(G_2)$ is written as
\begin{align*}
B_2=\left(
  \begin{array}{ccc}
  s-1 & n-2s+1 & s-1\\
  s & n-2s & 0\\
  s & 0 & 0\\
  \end{array}
\right),
\end{align*}
and the characteristic polynomial of $B_2$ is
\begin{align*}
\varphi_{B_2}(x)=x^{3}-(n-s-1)x^{2}-(n+s^{2}-2s)x+s^{2}n-sn-2s^{3}+2s^{2}.
\end{align*}
Note that the partition $V(G_2)=V(K_s)\cup V(K_{n-2s+1})\cup V((s-1)K_1)$ is equitable. In terms of Lemma 2.5, the largest root of $\varphi_{B_2}(x)=0$ equals $\rho(G_2)$.

Let $G_*=K_{\delta}\vee(K_{n-2\delta+1}\cup(\delta-1)K_1)$. We denote by $B_*$ the quotient matrix of $A(G_*)$ with respect to the partition $V(G_*)=V(K_{\delta})\cup V(K_{n-2\delta+1})\cup V((\delta-1)K_1)$.
Then
\begin{align*}
B_*=\left(
  \begin{array}{ccc}
  \delta-1 & n-2\delta+1 & \delta-1\\
  \delta & n-2\delta & 0\\
  \delta & 0 & 0\\
  \end{array}
\right),
\end{align*}
and its characteristic polynomial equals
\begin{align}\label{eq:4.3}
\varphi_{B_*}(x)=x^{3}-(n-\delta-1)x^{2}-(n+\delta^{2}-2\delta)x+\delta^{2}n-\delta n-2\delta^{3}+2\delta^{2}.
\end{align}
It is obvious that the partition $V(G_*)=V(K_{\delta})\cup V(K_{n-2\delta+1})\cup V((\delta-1)K_1)$ is equitable. From Lemma 2.5, the largest root of $\varphi_{B_*}(x)=0$ equals $\rho(G_*)$.

By a simple calculation, we obtain
\begin{align}\label{eq:4.4}
\varphi_{B_2}(x)-\varphi_{B_*}(x)=(s-\delta)f(x),
\end{align}
where $f(x)=x^{2}-(s+\delta-2)x+sn+\delta n-n-2s^{2}-2\delta s+2s-2\delta^{2}+2\delta$. Then the symmetry axis of $f(x)$ is $x=\frac{s+\delta-2}{2}$. According to $s\geq\delta+1$ and $n\geq2s$, we deduce
$\frac{s+\delta-2}{2}<s-1<n-\delta$. Hence, $f(x)$ is increasing with respect to $x\geq n-\delta$. For $x\geq n-\delta$, we possess
\begin{align}\label{eq:4.5}
f(x)=&x^{2}-(s+\delta-2)x+sn+\delta n-n-2s^{2}-2\delta s+2s-2\delta^{2}+2\delta\nonumber\\
\geq&(n-\delta)^{2}-(s+\delta-2)(n-\delta)+sn+\delta n-n-2s^{2}-2\delta s+2s-2\delta^{2}+2\delta\nonumber\\
=&-2s^{2}-(\delta-2)s+n^{2}-(2\delta-1)n\nonumber\\
\geq&-2\left(\frac{n}{2}\right)^{2}-(\delta-2)\left(\frac{n}{2}\right)+n^{2}-(2\delta-1)n \ \ \ \ \ \left(\mbox{since} \ s\leq\frac{n}{2} \ \mbox{and} \ \delta\geq2\right)\nonumber\\
=&\frac{1}{2}n^{2}-\left(\frac{5}{2}\delta-2\right)n\nonumber\\
\geq&\frac{1}{2}(5\delta-3)^{2}-\left(\frac{5}{2}\delta-2\right)(5\delta-3) \ \ \ \ \ (\mbox{since} \ n\geq5\delta-3)\nonumber\\
=&\frac{5\delta-3}{2}\nonumber\\
>&0 \ \ \ \ \ (\mbox{since} \ \delta\geq2).
\end{align}
It follows from \eqref{eq:4.4}, \eqref{eq:4.5} and $s\geq\delta+1$ that $\varphi_{B_2}(x)>\varphi_{B_*}(x)$ for $x\geq n-\delta$. Notice that $K_{n-\delta+1}$ is a proper subgraph of
$G_*=K_{\delta}\vee(K_{n-2\delta+1}\cup(\delta-1)K_1)$. Together with Lemma 2.3, we obtain $\rho(G_*)>\rho(K_{n-\delta+1})=n-\delta$, and so $\rho(G_2)<\rho(G_*)$. Together with \eqref{eq:4.1} and
\eqref{eq:4.2}, we conclude
$$
\rho(G)\leq\rho(G_1)\leq\rho(G_2)<\rho(G_*)=\rho(K_{\delta}\vee(K_{n-2\delta+1}\cup(\delta-1)K_1)),
$$
which contradicts $\rho(G)\geq\rho(K_{\delta}\vee(K_{n-2\delta+1}\cup(\delta-1)K_1))$.

\noindent{\bf Case 2.} $s=\delta$.

Recall that $G_1=K_s\vee(K_{n_1}\cup K_{n_2}\cup\cdots\cup K_{n_s})$. By Lemma 2.4 and $s=\delta$, we get
$$
\rho(G_1)\leq\rho(K_{\delta}\vee(K_{n-2\delta+1}\cup(\delta-1)K_1)),
$$
with equality holding if and only if $G_1=K_{\delta}\vee(K_{n-2\delta+1}\cup(\delta-1)K_1)$. Combining this with \eqref{eq:4.1}, we deduce
$$
\rho(G)\leq\rho(K_{\delta}\vee(K_{n-2\delta+1}\cup(\delta-1)K_1)),
$$
with equality if and only if $G=K_{\delta}\vee(K_{n-2\delta+1}\cup(\delta-1)K_1)$, a contradiction.

\noindent{\bf Case 3.} $s\leq\delta-1$.

Let $G_3=K_s\vee(K_{n-s-(\delta+1-s)(s-1)}\cup(s-1)K_{\delta+1-s})$. Recall that $G$ is a spanning subgraph of $G_1=K_s\vee(K_{n_1}\cup K_{n_2}\cup\cdots\cup K_{n_s})$, where $n_1\geq n_2\geq\cdots\geq n_s$
and $\sum\limits_{i=1}^{s}n_i=n-s$. Obviously, $n_s\geq\delta+1-s$ because $\delta(G_1)\geq\delta(G)=\delta$. In terms of Lemma 2.4, we conclude
\begin{align}\label{eq:4.6}
\rho(G_1)\leq\rho(G_3),
\end{align}
where the equality holds if and only if $(n_1,n_2,\ldots,n_s)=(n-s-(\delta+1-s)(s-1),\delta+1-s,\ldots,\delta+1-s)$. The quotient matrix of $A(G_3)$ according to the partition
$V(G_3)=V(K_s)\cup V(K_{n-s-(\delta+1-s)(s-1)})\cup V((s-1)K_{\delta+1-s})$ is equal to
\begin{align*}
B_3=\left(
  \begin{array}{ccc}
  s-1 & n-s-(\delta+1-s)(s-1) & (s-1)(\delta+1-s)\\
  s & n-s-(\delta+1-s)(s-1)-1 & 0\\
  s & 0 & \delta-s\\
  \end{array}
\right).
\end{align*}
The characteristic polynomial of $B_3$ is
\begin{align*}
\varphi_{B_3}(x)=&x^{3}-(n+s^{2}-\delta s-3s+2\delta-1)x^{2}\\
&+(\delta n-sn-n+\delta s^{2}-s^{2}-\delta^{2}s-\delta s+4s+\delta^{2}-2\delta)x\\
&+(\delta-s)(n-(\delta+1-s)(s-1)-1+s(s-1)(n-s-(\delta+1-s)(s-1)-1))\\
&+s(s-1)(n-s-(\delta+1-s)(s-1)-1).
\end{align*}
Since the partition $V(G_3)=V(K_s)\cup V(K_{n-s-(\delta+1-s)(s-1)})\cup V((s-1)K_{\delta+1-s})$ is equitable, it follows from Lemma 2.5 that $\rho(G_3)$ is the largest root of $\varphi_{B_3}(x)=0$.

Recall that $G_*=K_{\delta}\vee(K_{n-2\delta+1}\cup(\delta-1)K_1)$. Let $\theta=\rho(G_*)$. Recall that $\theta=\rho(G_*)>n-\delta$ and $\varphi_{B_*}(\theta)=0$. By plugging the value $\theta$ into $x$
of $\varphi_{B_3}(x)-\varphi_{B_*}(x)$, we conclude
\begin{align}\label{eq:4.7}
\varphi_{B_3}(\theta)=\varphi_{B_3}(\theta)-\varphi_{B_*}(\theta)=(\delta-s)g(\theta),
\end{align}
where $g(\theta)=(s-3)\theta^{2}+(n-\delta s+s+2\delta-4)\theta+s^{4}-(\delta+5)s^{3}+(n+2\delta+8)s^{2}-(2n+5)s+(2-\delta)n+2\delta^{2}-\delta$. Notice that
$$
-\frac{n-\delta s+s+2\delta-4}{2(s-3)}<n-\delta<\theta
$$
by $3\leq s\leq\delta-1$ and $n\geq\frac{1}{3}\delta^{2}+\delta$. This implies that
\begin{align}\label{eq:4.8}
g(\theta)>&g(n-\delta)\nonumber\\
=&(s-3)(n-\delta)^{2}+(n-\delta s+s+2\delta-4)(n-\delta)+s^{4}-(\delta+5)s^{3}\nonumber\\
&+(n+2\delta+8)s^{2}-(2n+5)s+(2-\delta)n+2\delta^{2}-\delta\nonumber\\
=&(s-2)n^{2}+(s^{2}-(3\delta+1)s+6\delta-2)n+s^{4}-(\delta+5)s^{3}\nonumber\\
&+(2\delta+8)s^{2}+(2\delta^{2}-\delta-5)s-3\delta^{2}+3\delta.
\end{align}
Let $h(x)=(s-2)x^{2}+(s^{2}-(3\delta+1)s+6\delta-2)x+s^{4}-(\delta+5)s^{3}+(2\delta+8)s^{2}+(2\delta^{2}-\delta-5)s-3\delta^{2}+3\delta$. Then we obtain the derivative function of $h(x)$ as
$$
h'(x)=2(s-2)x+s^{2}-(3\delta+1)s+6\delta-2.
$$
For $x\geq\frac{3\delta-s-1}{2}$, we deduce
$$
h'(x)=2(s-2)\cdot\frac{3\delta-s-1}{2}+s^{2}-(3\delta+1)s+6\delta-2=0.
$$
This implies that $h(x)$ is increasing in the interval $[\frac{3\delta-s-1}{2},+\infty)$. Note that $\frac{3\delta-s-1}{2}\leq\frac{3\delta-4}{2}<\frac{1}{3}\delta^{2}+\delta\leq n$ by $s\geq3$ and
$\delta\geq s+1$. Thus, we infer
\begin{align}\label{eq:4.9}
h(n)\geq&h\left(\frac{1}{3}\delta^{2}+\delta\right)\nonumber\\
=&\frac{\delta}{9}(\delta((s-2)\delta^{2}-3(s-2)\delta+3s^{2}-3s+3)-9s^{3}+27s^{2}-18s+9)\nonumber\\
&+s(s^{3}-5s^{2}+8s-5)\nonumber\\
\geq&\frac{\delta}{9}(\delta((s-2)(s+1)^{2}-3(s-2)(s+1)+3s^{2}-3s+3)-9s^{3}+27s^{2}-18s+9)\nonumber\\
&+s \ \ \ \ \ (\mbox{since} \ s\geq3 \ \mbox{and} \ \delta\geq s+1)\nonumber\\
=&\frac{\delta}{9}(\delta(s^{3}-3s+7)-9s^{3}+27s^{2}-18s+9)+s\nonumber\\
\geq&\frac{\delta}{9}((s+1)(s^{3}-3s+7)-9s^{3}+27s^{2}-18s+9)+s \ \ \ \ \ (\mbox{since} \ s\geq3 \ \mbox{and} \ \delta\geq s+1)\nonumber\\
=&\frac{\delta}{9}(s^{4}-8s^{3}+24s^{2}-14s+16)+s\nonumber\\
>&6\delta+s \ \ \ \ \ (\mbox{since} \ s\geq3)\nonumber\\
>&0.
\end{align}
It follows from \eqref{eq:4.7}, \eqref{eq:4.8}, \eqref{eq:4.9} and $\delta\geq s+1$ that
\begin{align}\label{eq:4.10}
\varphi_{B_3}(\theta)=(\delta-s)g(\theta)>(\delta-s)h(n)>0.
\end{align}
In what follows, we give the derivative function of $\varphi_{B_3}(x)$ as
\begin{align*}
\varphi_{B_3}'(x)=&3x^{2}-2(n+s^{2}-\delta s-3s+2\delta-1)x\\
&+\delta n-sn-n+\delta s^{2}-s^{2}-\delta^{2}s-\delta s+4s+\delta^{2}-2\delta.
\end{align*}
Note that $n\geq(\delta+1-s)(s-1)+(\delta+1-s)+s=s(\delta+2-s)$ and
$$
-\frac{-2(n+s^{2}-\delta s-3s+2\delta-1)}{2\times3}=\frac{n+s^{2}-\delta s-3s+2\delta-1}{3}<n-\delta,
$$
and so $\varphi_{B_3}'(x)$ is increasing in the interval $[n-\delta,+\infty)$. For $x\geq n-\delta$, we possess
\begin{align}\label{eq:4.11}
\varphi_{B_3}'(x)\geq&\varphi_{B_3}'(n-\delta)\nonumber\\
=&3(n-\delta)^{2}-2(n+s^{2}-\delta s-3s+2\delta-1)(n-\delta)\nonumber\\
&+\delta n-sn-n+\delta s^{2}-s^{2}-\delta^{2}s-\delta s+4s+\delta^{2}-2\delta\nonumber\\
=&n^{2}+(-2s^{2}+2\delta s+5s-7\delta+1)n\nonumber\\
&+3\delta s^{2}-s^{2}-3\delta^{2}s-7\delta s+4s+8\delta^{2}-4\delta\nonumber\\
\geq&\left(\frac{1}{3}\delta^{2}+\delta\right)^{2}+(-2s^{2}+2\delta s+5s-7\delta+1)\left(\frac{1}{3}\delta^{2}+\delta\right)\nonumber\\
&+3\delta s^{2}-s^{2}-3\delta^{2}s-7\delta s+4s+8\delta^{2}-4\delta \ \ \ \ \ (\mbox{since} \ n\geq\frac{1}{3}\delta^{2}+\delta)\nonumber\\
=&\frac{1}{9}(\delta^{2}(\delta^{2}+(6s-15)\delta-6s^{2}+6s+21)+(9s^{2}-18s-27)\delta-9s^{2}+36s)\nonumber\\
\geq&\frac{1}{9}(\delta^{2}((s+1)^{2}+(6s-15)(s+1)-6s^{2}+6s+21)\nonumber\\
&+(9s^{2}-18s-27)\delta-9s^{2}+36s) \ \ \ \ \ (\mbox{since} \ \delta\geq s+1)\nonumber\\
=&\frac{1}{9}((s^{2}-s+7)\delta^{2}+(9s^{2}-18s-27)\delta-9s^{2}+36s)\nonumber\\
\geq&\frac{1}{9}((s^{2}-s+7)(s+1)^{2}+(9s^{2}-18s-27)(s+1)-9s^{2}+36s) \ \ \ \ \ (\mbox{since} \ \delta\geq s+1\geq4)\nonumber\\
\geq&\frac{1}{9}(16(s^{2}-s+7)-9s^{2}+36s) \ \ \ \ \ (\mbox{since} \ s\geq3)\nonumber\\
=&\frac{1}{9}(7s^{2}+20s+112)\nonumber\\
>&0.
\end{align}
For $s=2$, one can check that inequalities \eqref{eq:4.10} and \eqref{eq:4.11} also hold. According to \eqref{eq:4.10} and \eqref{eq:4.11}, we infer
\begin{align}\label{eq:4.12}
\rho(G_3)<\theta=\rho(G_*).
\end{align}
Recall that $G_*=K_{\delta}\vee(K_{n-2\delta+1}\cup(\delta-1)K_1)$. Combining this with \eqref{eq:4.1}, \eqref{eq:4.6} and \eqref{eq:4.12}, we conclude
$$
\rho(G)\leq\rho(G_1)\leq\rho(G_3)<\rho(K_{\delta}\vee(K_{n-2\delta+1}\cup(\delta-1)K_1)),
$$
which is a contradiction to $\rho(G)\geq\rho(K_{\delta}\vee(K_{n-2\delta+1}\cup(\delta-1)K_1))$. Theorem 1.2 is proved. \hfill $\Box$

\section*{Declaration of competing interest}

The authors declare that they have no known competing financial interests or personal relationships that could have appeared to influence the work reported in this paper.

\section*{Data availability}

No data was used for the research described in the article.

\section*{Acknowledgments}

This work was supported by the Natural Science Foundation of Jiangsu Province (Grant No. BK20241949). Project ZR2023MA078 supported by Shandong Provincial Natural Science Foundation.

\end{document}